\documentclass[11pt,a4paper]{article}
\usepackage{authblk}
\usepackage[T1]{fontenc}
\usepackage[utf8]{inputenc}
\usepackage{lmodern,microtype,amsmath,amssymb,amsthm,booktabs,array}
\usepackage{tikz}
\usetikzlibrary{arrows.meta}
\usepackage[margin=27mm]{geometry}
\usepackage[colorlinks=true,linkcolor=blue,citecolor=blue,urlcolor=blue]{hyperref}
\usepackage{enumitem}
\setlist{nosep}
\newtheorem{theorem}{Theorem}[section]
\newtheorem{lemma}[theorem]{Lemma}
\newtheorem{proposition}[theorem]
{Proposition}
\newtheorem{Claim}[theorem]{Claim}
\newcommand{\ex}{\operatorname{ex}}

\title{The maximum number of triangles in a graph without pentagons}
\author[1]{Qing Cui}
\author[2]{Ervin Győri}
\author[1]{Chenke Zhang}
\author[1]{Xiutao Zhu}

\date{}

\affil[1]{School of Mathematics, Nanjing University of Aeronautics and Astronautics}
\affil[2]{Alfr\'ed R\'enyi Institute of Mathematics, Hungarian Academy of Sciences. }

\hypersetup{pdftitle={Triangles in a graph without pentagons}}
\begin{document}
\maketitle

\begin{abstract}
The generalized Tur\'an number $\ex(n,T,H)$ is the maximum number of $T$ in an $n$-vertex graph containing no $H$. Inspired by a conjecture of Erd\H{o}s on the function $\ex(n,C_5,K_3)$,  Bollob\'as and Gy\H ori proposed the problem of $\ex(n,K_3,C_5)$. In this paper, we prove that
$\ex(n,K_3,C_5)=(1/(3\sqrt3)+o(1))n^{3/2}$, finally solved this problem. 

As 
applications, we determine the asymptotic value of $\ex(n,\{\text{ind-}C_4,C_{5}\})$ which answers a question of Ergemlidze, Gy\H{o}ri and Methuku, and the  asymptotic value of  3-uniform hypergraph without Berge short cycle which solve a conjecture of Gy\H{o}ri and Salia.
\end{abstract}

\section{Introduction}

Tur\'an type problems are fundamental in extremal graph theory. A graph $G$ is said to be $F$-free for some graph $F$ if it does not contain $F$ as a subgraph. The ordinary Tur\'an number, $\ex(n,F)$  is the maximum number of edges an $F$-free simple graph on $n$ vertices can have. This has been of great interest ever since its start from the renowned theorem of Tur\'an \cite{T} in 1941, who determined $\ex(n,K_r)$, where $K_r$ is a complete graph on $r$ vertices. A natural generalization is the maximum  number of copies of another subgraph $H$ instead of edges while forbidding $F$, which is denoted by $\ex(n,H,F)$ and called generalized Tur\'an  number. Before initiating its general and systematic study by Alon and Shikhelman \cite{AS} in 2016, there had been sporadic results regarding this function, see \cite{GP}.

Although this concept was only introduced in 2016, many scholars had already studied this class of problems before that. For instance, Erd\H{o}s \cite{Erd62}, determined $\ex(n,K_r,K_t)$ for all $r<t$, and Gy\H{o}ri, Pach and Simonovits \cite{GyPS} studied $\ex(n,H,K_r)$, especially when $r=3$ and $H$ is bipartite. Probably the most interesting problems among all of these ones are those concerning pentagons $C_5$ (cycles of length 5) and triangles $C_3$. In 1984, Erd\H{o}s \cite{Erd84}  conjectured that the maximum number of $C_5$'s in a triangle free graph on $n$ vertices is at most $(n/5)^5$, which took a quarter of century until it was proved in 2013, see \cite{Grz, HHKNR}. Inspired by this conjecture, Bollob\'as and Gy\H{o}ri \cite{BG} studied the converse of this problem, $\ex(n,C_3,C_5)$. In their paper, they give a construction: take a $C_4$-free bipartite graph with each part has size $n/3$ first, and then double each vertex in one of the part and add an edge joining the old and the new copy. Using this graph, they showed the lower bound
$$\ex(n,C_3,C_5)\ge (\frac1{3\sqrt3}+o(1))n^{3/2}.$$ 
And then, several papers, see \cite{AS, CFSZ,EM, EMSG, FO} were written on this particular problem to find the matching upper bound. The best known result on this problem is as follows.

\begin{theorem}[Lv, He and Lu \cite{LHL}]\label{LHL}
$\ex(n,K_3,C_5)\le \frac{1}{2\sqrt6}n^{\frac{3}{2}}+o(n^{\frac{3}{2}})$.
\end{theorem}

After 19 years, we determine the iconic $\ex(n,K_3,C_5)$ in this paper.

\begin{theorem}\label{thm:main}
$\ex(n,K_3,C_5)=\left(\frac1{3\sqrt3}+o(1)\right)n^{3/2}.$
\end{theorem}

\subsection{Applications}
 We  also present two applications of our methods. 

A graph is
called induced-$F$-free if it does not contain an induced copy of $F$. In \cite{LTTZ}. Loh, Tait, Timmons, and Zhou introduced the problem of
 $\ex(n,\{\text{ind-}F,H\})$ and get an upper bound of $\ex(n,\{\text{ind-}K_{2,t},C_{2k+1}\})$. Later, Ergemlidze, Gy\H{o}ri and Methuku \cite{EGM} determined the asymptotically value of $\ex(n,\{\text{ind-}K_{2,t},C_{2k+1}\})$ except the case  $\ex(n,\{\text{ind-}C_4,C_{5}\})$. Hence they leave this case as an open problem.  In \cite{EM},  the authors showed $\ex(n,\{\text{ind-}C_4,C_{5}\})\le \frac{1}{2\sqrt[10]{2}}(1+o(1))n^{3/2}$, ant it was improved to $\ex(n,\{\text{ind-}C_4,C_{5}\})\le \frac{1}{\sqrt{6}}(1+o(1))n^{3/2}$ in \cite{LHL}. Here, using our method we solve this
problem completely.
\begin{theorem}\label{thm:induced}
    $\ex(n,\{\text{ind-}C_4,C_{5}\})= \frac{2}{3\sqrt{3}}(1+o(1))n^{3/2}$.
\end{theorem}

 One another application is about the 3-uniform hypergraph. 
A \emph{$3$-uniform hypergraph} $\mathcal{H}$ is a family of distinct
three-element subsets of a vertex set.
A \emph{Berge cycle} $\mathcal{B}C_\ell$ consists of distinct vertices
$v_1,\ldots,v_\ell$ and distinct hyperedges $h_1,\ldots,h_\ell$ such that
$\{v_i,v_{i+1}\}\subseteq h_i$, with indices modulo $\ell$.
For a family $\mathcal{F}$ of $3$-uniform hypergraphs, let
$\ex_3(n,\mathcal{F})$ denote the maximum number of
hyperedges in an $n$-vertex $\mathcal{F}$-free $3$-uniform hypergraph. In \cite{BG}, Bollob\'as and Gy\H{o}ri obtain an upper bound for $\ex_3(n,\mathcal{B}C_5)$. Later, Ergemlidze and Methuku \cite{EM} proved that, if $\mathcal{H}$ is a linear 3-uniform hypergraph and contains no $\{\mathcal{B}C_3,\mathcal{B}C_4,\mathcal{B}C_5\}$, then $e(\mathcal{H})\le \frac{1+o(1)}{6\sqrt{3}}n^{3/2}$.  Gy\H{o}ri and Salia \cite{GS} then conjecture that $e(\mathcal{H})=\frac{1+o(1)}{3\sqrt{3}}n^{3/2}$ without the linearity constraint. Using our result, we solve this conjecture.
\begin{theorem}\label{thm:berge}
Let  $\mathcal{B}=\{\mathcal{B}C_3,\mathcal{B}C_4,\mathcal{B}C_5\}$. Then $\ex_3(n,\mathcal{B})=\frac{1+o(1)}{3\sqrt{3}}n^{3/2}$.
\end{theorem}

The structure of the paper is as follows. In Section~\ref{sec:preliminaries}, we will give some preliminaries and  reduces the triangle problem to a weighted directed graph.
 In Section~\ref{sec:weighted}, we  will uses a partition argument and graph removal to get a new weighted graph whose total weights are still large.
In Section~\ref{sec:bipartite}, we prove the total weights can not be large in
bipartite graph. Then in Section~\ref{sec:completion}, we uses this bound to complete the proof of Theorem~\ref{thm:main}.
Section~\ref{sec:applications} proves Theorems~\ref{thm:induced}
and~\ref{thm:berge}.

\section{Preliminaries and a reduction}\label{sec:preliminaries}

Let $V$ be a set of size $n$. For a real function $h$ on $V^2$, define its 2-norm and cut norm by
\[\|h\|_2^2=n^{-2}\sum_{x,y}h(x,y)^2,\quad \|h\|_\square=\frac{1}{n^2}\max_{S,T\subseteq V}
 \left|\sum_{x\in S,\,y\in T}h(x,y)\right|.\]
For another function $W$ on $V^2$, let
\[
 (h\circ W)(x,y)=\frac1n\sum_z h(x,z)W(z,y),
 \qquad T_W=W\circ W.
\]
We will use some results from \cite{{CFSZ}}.
\begin{lemma}[Conlon, Fox, Sudakov and Zhao{\cite[Lemma~3.3]{CFSZ}}]\label{lem:cut-convolution}
Let $h:V^2\to\mathbb R$ and $W:V^2\to[0,\infty)$ satisfy
\[
 \frac1n\sum_y W(z,y)\le M\qquad(z\in V).
\]
Then $\|h\circ W\|_\square\le M\|h\|_\square$.
If $\frac1n\sum_y W(y,z)\le M$ also holds, then
$\|W\circ h\|_\square\le M\|h\|_\square$. 
\end{lemma}

\begin{lemma}[Conlon, Fox, Sudakov and Zhao~{\cite[Corollary~1.4]{CFSZ}}]\label{CFSZ}
Every $n$-vertex $C_5$-free graph can be made triangle-free by deleting $o(n^{3/2})$ edges.
\end{lemma}

For a symmetric function $W\ge0$ on $V^2$, write
\[
 t(C_\ell,W)=\frac1{n^\ell}\sum_{x_1,\ldots,x_\ell\in V}
             \prod_{i=1}^{\ell}W(x_i,x_{i+1}),
 \qquad x_{\ell+1}=x_1.
\]

\begin{theorem}[Conlon, Fox, Sudakov and Zhao~{\cite[Theorem~3.1]{CFSZ}}]\label{CFSZ3}
Fix $0<\varepsilon<1$ and $C\ge1$, and let $V$ be a probability
space. Suppose that $f\colon V^2\to[0,\infty)$ and
$ g\colon V^2\to[0,1]$ are measurable symmetric functions with
\[
 t(C_4,f)\le C
 \qquad\text{and}\qquad
 \|f-g\|_\square\le\varepsilon^4.
\]
Then
\[
t(C_5,g)\le  t(C_5,f)+11C\varepsilon.
\]
\end{theorem}

\begin{theorem}[Conlon, Fox, Sudakov and Zhao~{\cite[Theorem~4.1]{CFSZ}}]
\label{CFSZ2}
For each $\zeta>0$, there exists $\delta>0$ with the following property.
If $g:V^2\to[0,1]$ is symmetric and $t(C_5,g)\le\delta$, then
there is a symmetric set $\mathcal D\subseteq V^2$ such that
\[
 t(C_5,g\mathbf1_{\mathcal D^c})=0,
 \qquad \frac1{n^2}\sum_{(x,y)\in\mathcal D}g(x,y)\le\zeta.
\]
\end{theorem}

Let $k:V^2\to[0,\infty)$ be a function and write $K(x,y)=k(x,y)+k(y,x)$. The support of $k$ is the graph on $V$ such that $xy$ is an edge if and only if $K(x,y)>0$, including a loop if $x=y$.  In a directed graph $D$,  we assign a coefficient to each $2$-walk $izj$ in its support: the coefficient is $1/2$ if its orientation is $i\leftarrow z\rightarrow j$, and $1$ otherwise. For $x,y\in V$, define the $2$-walk capacity of $(x,y)$ by
\begin{equation*}
\begin{split}
 \mathcal{C}_k(x,y)=\frac{1}{n}\sum_z\bigl(&k(x,z)k(y,z)+\tfrac12k(z,x)k(z,y)\\
                     &+k(x,z)k(z,y)+k(z,x)k(y,z)\bigr).
\end{split}
\end{equation*}
Several direct consequences of the capacity bound will be used. 
\begin{lemma}
Let $K(x,y)=k(x,y)+k(y,x)$. If $\mathcal{C}_k(x,y)\le 1$, then  
\begin{equation}\label{eq:capacitybounds}
\|K\|_2^2\le\frac43,\qquad \frac1n\sum_z K(x,z)^2\le2.
\end{equation}

\end{lemma}
 
 \begin{proof}
 When summing over $x$ of $ 
\mathcal{C}_k(x,x)= \frac1n\sum_z\bigl(k(x,z)^2+\tfrac12k(z,x)^2+2k(x,z)k(z,x)\bigr)\le 1$, we get

$$\tfrac32\|k\|_2^2+\frac2{n^2}\sum_{x,y}k(x,y)k(y,x)\le1,$$ 
which implies
\begin{equation*}
\|K\|_2^2=2\|k\|_2^2+\frac2{n^2}\sum_{x,y}k(x,y)k(y,x)\le \frac{4}{3}.
\end{equation*}
On the other hand, a direct calculation shows $ \frac{1}{n}\sum_z K(x,z)^2\le 2\mathcal{C}_k(x,x)\le 2$  We are done.    
 \end{proof}

Next, we will reduce the problem to a weighted graph with $\{K_3,C_5\}$-free support and study this weighted graph to get contradiction. Two triangles are in the same \emph{triangle block} if they can be joined by a sequence of triangles in which consecutive triangles share an edge. A \emph{book} has triangles $abz_1,\ldots,abz_s$ with common spine $ab$ and page vertices $z_1,\ldots,z_s$.

\begin{proposition}[\cite{EMSG,EM}]\label{prop:blocks}
In a $C_5$-free graph in which every edge lies in a triangle, every triangle block is a book or a $K_4$. Any two nonadjacent vertices have at most two common neighbors; if there are two, those neighbors are adjacent.
\end{proposition}

\begin{lemma}\label{deletion}
 If the upper bound in Theorem~\ref{thm:main} fails, then for arbitrarily large $n$ we can find a $\{K_3,C_5\}$-free graph $H$ of order $n$ with maximum degree $\Delta(H)\le24\sqrt{n}$.   
\end{lemma}
\begin{proof}
If the upper bound in Theorem~\ref{thm:main} fails, there exists a fixed $\eta>0$ such that we can find arbitrarily large $C_5$-free graphs $G_N$ of order $N$ with $t(G_N)\ge(\frac{1}{3\sqrt{3}}+\eta)N^{3/2}$, here $t(G)$ denotes the number of triangles in $G$.
First remove the edges of $G_N$ in no triangle. At each intermediate step $G_i$, if there exists a vertex which belonging to less  than $\frac{1}{3\sqrt{3}}(i^{3/2}-(i-1)^{3/2})$ triangles, then delete this vertex and remove edges in no remaining triangles, say the obtained graph being $G_{i-1}$. The process terminates at the graph $G_{n}$. 

Since at each step,  $t(G_{i-1})\ge t(G_i)-\frac{1}{3\sqrt{3}}(i^{3/2}-(i-1)^{3/2})$ and $t(G_N)\ge(\frac{1}{3\sqrt{3}}+\eta)N^{3/2}$, we have
$t(G_{n})\ge (\frac{1}{3\sqrt{3}})n^{3/2}+\eta N^{3/2}\ge (\frac{1}{3\sqrt{3}}+\eta)n^{3/2}$. Also by Theorem \ref{LHL}, there is a big constant $C$ such that $t(G_{n})\le Cn^{3/2}$. This implies that $(C-\frac{1}{3\sqrt{3}})n^{3/2}\ge \eta N^{3/2}$ and so $n$ is also large enough.

For every vertex $x$, write $N(x)=N_{G_n}(x)$ and $d(x)=d_{G_n}(x)$.
The graph $G_n[N(x)]$ induced by $N(x)$ contains no path of length three, since such a path together with $x$ would form a $C_5$. Its components are stars or triangles, so the number of triangles containing $x$ is $e(G_n[N(x)])\le d(x)$.
At the final step, $\delta(G_n)\ge \frac{1}{3\sqrt{3}}(n^{3/2}-(n-1)^{3/2})\ge \frac{1}{3\sqrt{3}}\sqrt{n}$. For any vertex $x$,  let $N_2(x)$ be the set of vertices  at distance two from $x$. Every vertex in $N_2(x)$ has at most two neighbors in $N(x)$ by Proposition~\ref{prop:blocks}. Then
\[
2n\ge  e(N(x),N_2(x))=\sum_{v\in N(x)}d(v)-d(x)-2e(G_n[N(x)])\ge d(x)(\delta(G_n)-3).
\]
 Since $n$ is large, we have $\delta(G_n)-3\ge\sqrt n/(6\sqrt3)$, so $d_{G_n}(x)\le12\sqrt3\,\sqrt n<24\sqrt n$.

When we choose one triangle from each $K_4$-block, these triangles are edge-disjoint. Thus by Proposition \ref{CFSZ},  there are $o(n^{3/2})$ $K_4$-blocks. Deleting the edges in each such block removes exactly four triangles. Therefore at least $(\frac{1}{3\sqrt{3}}+\eta/2)n^{3/2}$ triangles remain for sufficiently large $n$, and all remaining blocks are books. Then we delete all spines in these book-blocks, the resulting graph, say $H$, hence is $\{K_3, C_5\}$-free graph.
\end{proof}

\begin{figure}[htbp]
\centering
\begin{tikzpicture}[scale=.85,every node/.style={font=\small},vertex/.style={circle,fill,inner sep=1.8pt}]
 \node[vertex,label=left:$a$] (a) at (0,1) {};
 \node[vertex,label=left:$b$] (b) at (0,-1) {};
 \foreach \i/\y in {1/1.5,2/0,3/-1.5}{
  \node[vertex,label=right:$z_\i$] (z\i) at (2,\y) {};
  \draw (a)--(z\i)--(b);
 }
 \draw[thick] (a)--(b);
 \node at (1,-2.1) {a book};
 \draw[-{Stealth[length=2mm]}] (3.2,0)--node[above]{delete $ab$}(5.1,0);
 \begin{scope}[xshift=6cm]
  \node[vertex,label=left:$a$] (aa) at (0,1) {};
  \node[vertex,label=left:$b$] (bb) at (0,-1) {};
  \foreach \i/\y in {1/1.5,2/0,3/-1.5}{
   \node[vertex,label=right:$z_\i$] (zz\i) at (2,\y) {};
   \draw[-{Stealth[length=1.6mm]}] (aa)--(zz\i);
   \draw[-{Stealth[length=1.6mm]}] (bb)--(zz\i);
  }
  \node at (1,-2.1) {its oriented sides};
 \end{scope}
\end{tikzpicture}
\caption{Each triangle contributes two arcs.}\label{fig:book}
\end{figure}
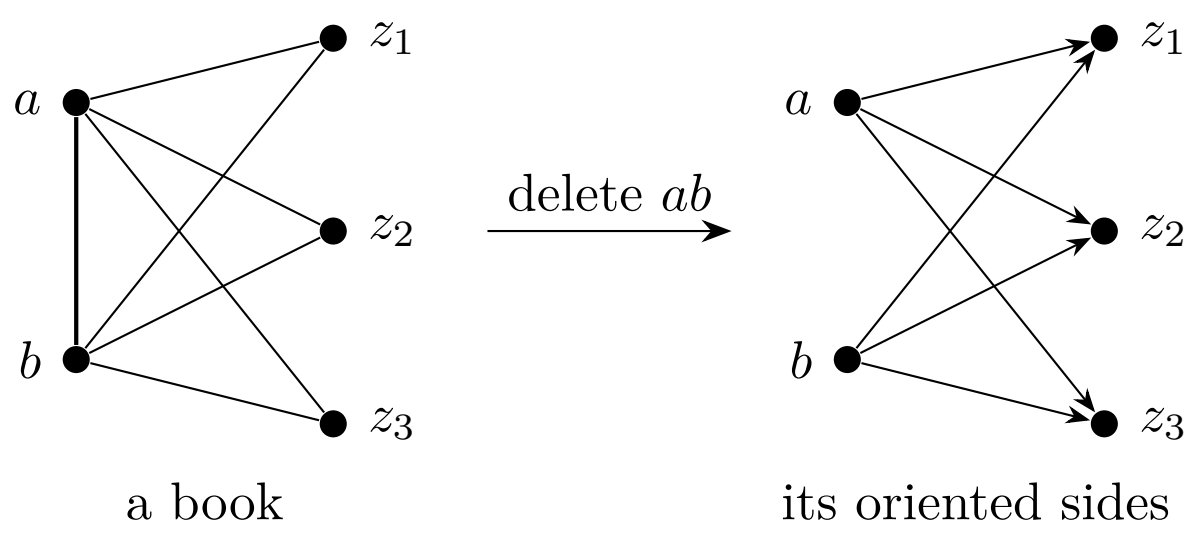

Let $H$ be the graph in the preceding proof, with the same constant $\eta>0$. For every remaining book-block with spine $ab$ and page vertices $z_1,\ldots,z_s$, orient the two sides $az_i,bz_i$ from the spine endpoints to the page vertex, as in Figure~\ref{fig:book}. This gives an oriented graph $D$ whose underlying graph is $H$. Write $V$ for its vertex set and define $f:V^2\to\{0,\sqrt n\}$ by $f(x,y)=\sqrt n$ if $xy$ is an arc in $D$ and $f(x,y)=0$ otherwise. Then we have the following properties. 

\begin{lemma}\label{lem:bookreduction} There is a nonnegative function $R$ such that
\[
 \frac1{n^2}\sum_{x,y}f(x,y)\ge\frac{2}{3\sqrt{3}}+\eta ,\qquad
 \mathcal{C}_f(i,j)\le1+R(i,j),\qquad
 \max_i\frac{1}{n}\sum_jR(i,j)\le\frac{48}{\sqrt{n}}.
\]
\end{lemma}
\begin{proof}
Let $G$ be the graph after the $K_4$-blocks have been removed and before the spines are deleted. Each of its triangles contributes two arcs of weight $\sqrt n$ to $D$. Therefore
\[
 \frac1{n^2}\sum_{x,y}f(x,y)=\frac{|E(D)|}{n^{3/2}}
 =\frac{2\tau(G)}{n^{3/2}}\ge\frac{2}{3\sqrt{3}}+\eta.
\]
For any two nonadjacent vertices $i,j$ in $G$,  they have at most two common neighbors by Proposition \ref{prop:blocks}.  If there is at most one common neighbor, then $\mathcal{C}_f(i,j)$ is at most one. If there are two common neighbors, then the two neighbors form the spine of a book block by Proposition \ref{prop:blocks}. The two $2$-walks both have orientation $i\leftarrow z \rightarrow j$, and hence each has coefficient $1/2$. For an adjacent pair which is not a spine in $G$, they have no common neighbor in the directed graph $D$. So in these three cases, we can let $R(i,j)=R(j,i)=0$.  For a spine $ij$ with $s$ pages, define $R(i,j)=R(j,i)=s$. In the case $i=j$, define $R(i,i)=d_D^+(i)+d_D^-(i)/2$.
Here $d_D^+(i)$ and $d_D^-(i)$ denote the outdegree and indegree of $i$, respectively

Above definition satisfies $\mathcal{C}_f(i,j)\le 1+R(i,j)$. Also for every vertex $i$, 
\[\sum_jR(i,j)=\sum_{ij ~\text{is spine}}R(i,j)+R(i,i)=d_D^+(i)+d_D^+(i)+\frac{d_D^-(i)}{2}.\]
Consequently, $\sum_j R(i,j)\le2\Delta(G_n)\le 48\sqrt{n}$ and we are done.
\end{proof}

\section{A weighted reduction for book graphs}\label{sec:weighted}
Take $D$, $V$, $R$, $\eta$ and $f$ from Section~2, and write
$n=|V|$ and $r=48/\sqrt{n}$.

\begin{proposition}\label{prop:reduction}
For sufficiently large $n$ in terms of $\eta$, there is a function
$k\ge0$ on $V^2$ such that its support has no closed $5$-walk and 
\[
 \frac1{n^2}\sum_{x,y}k(x,y)>\frac{2}{3\sqrt{3}},
 \qquad \mathcal{C}_k(x,y)\le1\quad(x,y\in V).
\]
\end{proposition}

\begin{proof}
Let $\mathcal{P}=\{P_1,\ldots,P_m\}$ be any partition of $V$ and
$p_i=|P_i|/n$. Its block average is the function
\[
 g(x,y)=g_{ij}=\frac1{|P_i||P_j|}\sum_{u\in P_i,\,v\in P_j}f(u,v)
 \qquad(x\in P_i,\ y\in P_j).
\]
First we prove two inequalities which will be used many times. By Cauchy--Schwarz on each $P_j$ and Lemma~\ref{lem:bookreduction}, we have
\begin{equation}\label{2-norm bound}
  \begin{aligned}
\|g\|_2^2
 &\le\sum_i\frac{p_i}{n}\sum_z
          \left(\frac1{|P_i|}\sum_{x\in P_i}f(x,z)\right)^2=\sum_i\frac{p_i}{|P_i|^2}\sum_{x,y\in P_i}
          \frac1n\sum_z f(x,z)f(y,z)\\
 &\le\sum_i\frac{p_i}{|P_i|^2}\sum_{x,y\in P_i}\bigl(1+R(x,y)\bigr)
 \le\sum_i(p_i+r)=1+mr.
\end{aligned}  
\end{equation}
The last inequality uses $\sum_{x,y\in P_i}R(x,y)\le |P_i|nr$. Moreover, if $\mathcal{P}'$ is a refinement of $\mathcal{P}$, then for any two parts $P_i,P_j\in \mathcal{P}$, we have $\sum_{x\in P_i,\,y\in P_j}(g'-g)(x,y)=0$ and so
\[
 \sum_{x,y}g(x,y)(g'-g)(x,y)
 =\sum_{i,j}g_{ij}\sum_{x\in P_i,\,y\in P_j}(g'-g)(x,y)=0.
\]
Therefore
\begin{equation}\label{2-norm equality}
    \|g'\|_2^2- \|g\|_2^2=\|g'-g\|_2^2+\frac2{n^2}\sum_{x,y}g^{i,j}(x,y)(g'-g)(x,y)
 =\|g'-g\|_2^2.
\end{equation}
\begin{Claim}\label{cl:partitions}
Fix $0<\epsilon<1/3$. For sufficiently large $n$, there are partitions
$\mathcal P$  and its refinement $\mathcal P'$, such that $m:=|\mathcal{P}|$ and $|
\mathcal{P'}|$ depend only on $\epsilon$.
Their block averages $g,g'$ satisfy
\begin{equation*}
 \|f-g'\|_\square\le\frac{\epsilon^3}{168m^2},
 \qquad \|g'\|_2^2-\|g\|_2^2\le\epsilon^2.
\end{equation*}
\end{Claim}

\begin{proof}
 We will prove this claim by finding a sequence of partitions $\mathcal{P}^0=\{V\}, \mathcal{P}^1,\ldots,\mathcal{P}^L$ such that they satisfy: $\mathcal{P}^{i+1}$ is a refinement of $\mathcal{P}^{i}$,  $m_i:=|\mathcal{P}^i|$ depends only on $\epsilon$, $\|f-g^i\|_\square\le \frac{\epsilon^3}{168m_{i-1}^2}$ and $L\le \lceil\frac{2}{\epsilon^2}\rceil$.  Then we will show in this sequences  $\mathcal{P}^{L-1}, \mathcal{P}^{L}$ are what we want.

First we start from $\mathcal{P}^0=\{V\}$ and suppose have found the sequence of partitions $\mathcal{P}^0, \mathcal{P}^1,\ldots,\mathcal{P}^{i-1}, \mathcal{P}^i=\{P_1,P_2,\ldots,P_{m_i}\}$. Let $g^{i}$ denote the block average function of this $\mathcal{P}^i$ and $\delta:=\frac{\epsilon^3}{168m_i^2}$. If $\|f-g^i\|_\square\le\delta$ already holds, then let $\mathcal{P}^{i+1}=\mathcal{P}^{i}$. If $\|f-g^i\|_\square>\delta$, then choose $S,T$ with $|n^{-2}\sum_{x\in S,\,y\in T}(f-g^i)(x,y)|>\delta$ and refine each $P_t$ in $\mathcal{P}^i$ by 
$$\{P_t\cap(S\cap T), P_t\cap(S\setminus T),P_t\cap(T\setminus S), P_t\setminus(S\cup T)\}.$$
Let this refinement be $\mathcal{P}^{i,1}$ and $g^{i,1}$ be its block average. We check the first condition. If $\|f-g^{i,1}\|_\square\le\delta$, then we define $\mathcal{P}^{i+1}=\mathcal{P}^{i,1}$. If $\|f-g^{i,1}\|_\square>\delta$, then find a new $S\times T$ such that $|n^{-2}\sum_{x\in S,\,y\in T}(f-g^{i,1})(x,y)|>\delta$ and refine each part in $\mathcal{P}^{i,1}$ by the same rule as above turn. The process continues until at some refinement $\mathcal{P}^{i,t}$ when $\|f-g^{i,t}\|_\square\le \delta$ and then define $\mathcal{P}^{i+1}=\mathcal{P}^{i,t}$. Next we prove $t\le \left\lceil2/\delta^2\right\rceil$. Suppose not, note that during the refinement process from $\mathcal{P}^{i,j}$ to $\mathcal{P}^{i,j+1}$,  
In particular, (\ref{2-norm equality})and Cauchy--Schwarz gives
\[
\begin{aligned}
 &\|g^{i,j+1}\|_2^2- \|g^{i,j}\|_2^2=\|g^{i,j+1}-g^{i,j}\|_2^2\\
 \ge&\left|\frac1{n^2}\sum_{x\in S,\,y\in T}(g^{i,j+1}-g^{i,j})(x,y)\right|^2=\left|\frac1{n^2}\sum_{x\in S,\,y\in T}(f-g^{i,j})(x,y)\right|^2
 >\delta ^2.
\end{aligned}
\]
Then at step $j=\lceil2/\delta^2\rceil$, $\|g^{i,j+1}\|_2^2\ge \delta^2\lceil2/\delta^2\rceil>2$. However,  since each time we refine $\mathcal{P}^{i,j}$, every part in $\mathcal{P}^{i,j}$ is splitted into at most four new parts, we have $|\mathcal P^{i,j+1}|\le m_{i}4^{j+1}$  and then (\ref{2-norm bound}) tells us  $\|g^{i,j+1}\|\le 1+|\mathcal{P}^{i,j+1}|r<2$ when $n$ is large enough, a contradiction. So $t\le \left\lceil2(168m_i^2/\epsilon^3)^2\right\rceil$, which also implies $|P^{i+1}|\le m_i\cdot4^{t}$, depends only on $\epsilon$.

Now we get a longer sequence $\mathcal{P}^0, \mathcal{P}^1,\ldots,\mathcal{P}^{i-1}, \mathcal{P}^{i+1}$ with $\|f-g^{i+1}\|_\square\le \frac{\epsilon^3}{168m_i^2}$. Next  we check the second condition. If $\|g^{i+1}\|_2^2-\|g^i\|_2^2\le \epsilon^2$, then the sequence ends and $\mathcal{P}^{i},\mathcal{P}^{i+1}$ are two partitions satisfy the claim.  If $\|g^{i+1}\|_2^2-\|g^i\|_2^2> \epsilon^2$, then we repeat the refinement on $\mathcal{P}^{i+1}$ as the last turn with $\delta=\epsilon^3/{168m_{i+1}^2}$.
 The process continues until we find $\mathcal{P}^0=\{V\}, \mathcal{P}^1,\ldots,\mathcal{P}^L$ with $\|g^L\|_2^2-\|g^{L-1}\|\le \epsilon^2$.
As definition of this sequence, all $\mathcal P^{i}$ are different except the possible case $\mathcal P^{L}=\mathcal P^{L-1}$. Since $\|g^{L}\|_2^2\ge L\epsilon^2$, the same reason shows $L\le \lceil2/\epsilon^2\rceil$, and we are done.
\end{proof}

Take the partitions $\mathcal{P}=\{P_1,P_2,\ldots, P_m\}$ and $\mathcal{P'}$ in Claim~\ref{cl:partitions}. Let $p_i:=|P_i|/n$ and $h:=g'-g$ and
\[
 e_i=\frac1{n|P_i|}\sum_{x\in P_i,\,z}
                 \bigl(h(x,z)^2+h(z,x)^2\bigr).
\]
Retain the parts with $p_i\ge\epsilon/m$ and $e_i\le\epsilon$, and
let $Z$ be the union of the other parts. Write
$F(x,y)=f(x,y)+f(y,x)$.

\begin{Claim}\label{cl:retained-capacity}
We have $|Z|\le3\epsilon n$. For $x,y$ in retained parts,
$\mathcal{C}_g(x,y)\le1+4\epsilon$.
\end{Claim}

\begin{proof}
By \eqref{2-norm equality} and Claim \ref{cl:partitions}, $\sum_i p_i e_i=2\|h\|_2^2=2(\|g'\|_2^2-\|g\|_2^2)\le2\epsilon^2$.
Hence
\[
 \frac{|Z|}{n}
 \le\sum_{p_i<\epsilon/m}p_i+\sum_{e_i>\epsilon}p_i
 \le\epsilon+\frac1\epsilon\sum_i p_i e_i
 \le3\epsilon.
\]

To show $\mathcal{C}_g(x,y)\le 1+4\epsilon$, recall the function $g'=g+h$.
For example, for every $P_\ell$,
\[
 \sum_{\substack{x\in P_i,\,y\in P_j\\ z\in P_\ell}}g(x,z)h(y,z)
 =|P_i|g_{i\ell}\sum_{v\in P_j,\,z\in P_\ell}h(v,z)=0.
\]
The same equalty holds for $g(x,z)h(z,y),g(z,x)h(y,z), g(z,x)h(z,y)$. Then for $x\in P_i,y\in P_j$
a direct calculation shows.
\begin{equation}\label{1}
 \frac1{|P_i||P_j|}\sum_{x\in P_i, y\in P_j}\mathcal{C}_{g'}(x,y)
 =\mathcal{C}_g(x,y)
 +\frac1{|P_i||P_j|}\sum_{x\in P_i,\,y\in P_j}\mathcal{C}_h(x,y).   
\end{equation}

Next consider the cut norm $\|\mathcal{C}_{g'}-\mathcal{C}_f\|_\square$. For convenient, let $f^{res}(x,y)=f(y,x)$ and $g'^{res}(x,y)=g'(y,x)$. Then
 \begin{equation*}
 \begin{split}
   \|\mathcal{C}_{g'}-\mathcal{C}_f\|_\square=\frac{1}{n^2}\max_{S,T}|\sum_{x\in S,y\in T}&(g'\circ g'-f\circ f)(x,y)+(g'^{res}\circ g'^{res}-f^{res}\circ f^{res})(x,y)\\
   &+(g'\circ g'^{res}-f\circ f^{res})(x,y)+\frac{1}{2}(g'^{res}\circ g'-f^{res}\circ f)(x,y)|.  
 \end{split}  
 \end{equation*}
Use the first term as an example. We know $n^{-1}\sum_y F(x,y)\le24$ by $\Delta(D)\le 24\sqrt{n}$ and for $x\in P_i$ we have, 
\begin{equation}\label{eq:block-row-bound}
\frac1n\sum_y\bigl(g(x,y)+g(y,x)\bigr)
 =\frac1{|P_i|}\sum_{u\in P_i}\frac1n\sum_yF(u,y)
 \le24.
\end{equation}
The same inequality holds for $g'$ and hence from Lemma~\ref{lem:cut-convolution}, the first term satisfies
\begin{equation*}
    \begin{split}
  \|g'\circ g'-f\circ f\|_\square&\le \|(g'-f)\circ g'\|_\square +\|f\circ (g'-f)\|_\square \\
    &\le 24\|(g'-f)\|_\square +24\|(g'-f)\|_\square\le 48\frac{\epsilon^3}{168m^2}.
    \end{split}
\end{equation*}
The other three terms in $\|\mathcal{C}_{g'}-\mathcal{C}_f\|_\square$ have similar bounds,  so we have 
 \[
 \|\mathcal{C}_{g'}-\mathcal{C}_f\|_\square
 \le (1+\tfrac12+1+1)\cdot48\frac{\epsilon^3}{168m^2}=\frac{\epsilon^3}{m^2}.
\]
Then in the  retained parts $P_i,P_j$,
\begin{equation}\label{2}
   \begin{aligned}
 &\frac1{|P_i||P_j|}\sum_{u\in P_i,\,v\in P_j}\mathcal{C}_{g'}(u,v)
= \frac1{|P_i||P_j|}\sum_{u\in P_i,\,v\in P_j}(\mathcal{C}_{f}+\mathcal{C}_{g'}-\mathcal{C}_{f})(u,v)\\
 \le &\frac{\sum_{u\in P_i,\,v\in P_j}(1+R(u,v))+ n^2\|\mathcal{C}_{g'}-\mathcal{C}_f\|_\square}{|P_i||P_j|}
 \le 1+\frac{r}{p_j}+\frac{\epsilon^3/m^2}{p_i p_j}
 \le 1+2\epsilon.
\end{aligned} 
\end{equation}
In these inequalities, we used Lemma \ref{lem:bookreduction}, $p_i,p_j\ge\epsilon/m$ and
$mr\le\epsilon^2$($n$ is large).

On the other hand, by Cauchy--Schwarz and $(a+b)^2\le2(a^2+b^2)$,
\[
\begin{aligned}
 &\left|\frac1{|P_i||P_j|}\sum_{u\in P_i,\,v\in P_j}\mathcal{C}_h(u,v)\right|\le\sum_{u\in P_i,\,v\in P_j,\,z}\frac{(|h(u,z)|+|h(z,u)|(|h(v,z)|+|h(z,v)|)}{n|P_i||P_j|}\\
&\quad\le\sqrt{2e_i}\sqrt{2e_j}\le2\epsilon.
\end{aligned}
\]
Therefore, combining this and (\ref{1}), (\ref{2}), we have $\mathcal{C}_g(x,y)\le1+4\epsilon$ on
the retained parts.
\end{proof}

Fix $T\ge1$ and write $G(x,y)=g(x,y)+g(y,x)$. For \(x\in P_i\) and \(y\in P_j\), define
\[
M(x,y)=
\begin{cases}
1, & \text{if both } P_i \text{ and } P_j \text{ are retained and }
    G(x,y)\le T,\\
0, & \text{otherwise}.
\end{cases}
\] 
and let 
 $k_1=\frac{M}{1+2\epsilon}g$, $K_1(x,y)=k_1(x,y)+k_1(y,x)$.
We have, 
\begin{Claim}
The function $k_1$ satisfies $\mathcal{C}_{k_1}\le1$.
Moreover,
\begin{equation}\label{eq:bounded-weights}
 \frac1{n^2}\sum_{x,y}(f-k_1)(x,y)\le96\epsilon+\frac4T.
\end{equation}
\end{Claim}

\begin{proof}
If $x$ or $y$ belongs to $Z$, then $\mathcal{C}_{k_1}(x,y)=0$.
Otherwise Claim~\ref{cl:retained-capacity} gives
\[
 \mathcal{C}_{k_1}(x,y)
 \le\frac{\mathcal{C}_g(x,y)}{(1+2\epsilon)^2}
 \le\frac{1+4\epsilon}{(1+2\epsilon)^2}\le1.
\]

Since $\sum_{x,y}f(x,y)=\sum_{x,y}g(x,y)$, the total losses of the weights satisfy
\[
\begin{aligned}
&\frac1{n^2}\sum_{x,y}(f-k_1)(x,y)=\frac1{n^2}\sum_{x,y}(g-k_1)(x,y)=\frac1{n^2}\sum_{x,y}\biggl(g(1-M)
        +\frac{2\epsilon}{1+2\epsilon}gM\biggr)(x,y)
\\
=&\frac{1}{n^2}\sum_{x\in Z,\,y\in V}G(x,y)+\frac{1}{2n^2}\sum_{G(x,y)>T}G(x,y)   +\frac{2\epsilon M}{1+2\epsilon}\frac{1}{n^2}\sum_{x,y}g(x,y).
\end{aligned}
\]
By \eqref{eq:block-row-bound}, we have $\frac{1}{n^2}\sum_{x\in Z,\,y\in V}G(x,y)\le 24\frac{|Z|}{n}$ and $\frac1{n^2}\sum_{x,y}g(x,y)
 =\frac1{2n^2}\sum_{x,y}G(x,y)\le 12$. Then the 2-norm bound in (\ref{2-norm bound}) gives, 
 \[\frac{1}{2n^2}\sum_{G(x,y)>T}G(x,y)\le\frac{1}{2Tn^2}\sum_{G(x,y)>T}G^2(x,y)\le \frac{\|G\|^2_2}{(2T)}\le \frac{4\|g\|_2^2}{2T}\le \frac{4}{T}.\]
Combining the bound of $Z$ in Claim \ref{cl:retained-capacity}, we get
\[
\frac1{n^2}\sum_{x,y}(f-k_1)(x,y)\le72\epsilon+\frac{4}{T}+2\epsilon\cdot12
 =96\epsilon+\frac{4}{T}.    
\]
The proof is done.\end{proof}
 Now the function $k_1$ already satisfies $\mathcal{C}_{k_1}\le 1$ and the total weights  do not loss too much compared to $f$. Next we continue to revise $k_1$ to make its  support containing no closed $5$-walk.
\begin{Claim}\label{cl:removing-walks}
For fixed $T$ and sufficiently small $\epsilon$, there is a function
$0\le k\le k_1$ whose support has no closed $5$-walk and
\[
 \frac1{n^2}\sum_{x,y}(k_1-k)(x,y)\le\frac{\eta}{4}.
\]
\end{Claim}

\begin{proof}
Let $\widetilde F=\frac{M}{1+2\epsilon}F$. As $R(x,y)\le nr$ and $n^{-2}\sum_{x,y}R(x,y)\le r$ in Lemma \ref{lem:bookreduction}, we have
\[
\|R\|_2^2\le(nr)\frac1{n^2}\sum_{x,y}R(x,y)\le nr^2=48^2.
\]
Expanding $T_F$, together with Lemma \ref{lem:bookreduction}, implies
$T_F(x,y)\le2+2R(x,y).$
Hence $\|T_F\|_2^2\le(2+2\|R\|_2)^2\le98^2$.
Using symmetry of $F$ and $0\le \widetilde F/T\le F$, we get
\[
\begin{aligned}
 t(C_4,\widetilde F/T)&\le t(C_4,F)
 =\frac1{n^2}\sum_{x,y}\left(\frac1n\sum_z F(x,z)F(z,y)\right)^2\\
 &=\|T_F\|_2^2\le(2+2\|R\|_2)^2\le98^2.
\end{aligned}
\]

For any $S_1,S_2\subseteq V$, the indicator $M$ is constant on
$P_i\times P_j$, so
\[
 \frac1{n^2}\left|\sum_{x\in S_1,\,y\in S_2}
                 \bigl(M(f-g')\bigr)(x,y)\right|
 \le\sum_{i,j}\frac1{n^2}
       \left|\sum_{x\in S_1\cap P_i,\,y\in S_2\cap P_j}
                      (f-g')(x,y)\right|
 \le \frac{\epsilon^3}{168}.
\]
The second inequality comes from Claim \ref{cl:partitions}. Also $\|Mh\|_\square\le n^{-2}\sum_{x,y}|h(x,y)|\le\|h\|_2
\le\epsilon$. Then we have 
\[
 \|\frac{\widetilde F}{T}-\frac{K_1}{T}\|_\square
 \le\frac{2}{T}\|M(f-g)\|_\square\le \frac{2}{T}\|M(f-g)\|_\square
 =\frac{2}{T}\left(\frac{\epsilon^3}{168}+\epsilon\right)\le\frac{3\epsilon}{T}.
\]

Also by the definition of $K_1(x,y)$, we know $\frac{K_1}{T}\le 1$. Above two results show that the function  $\frac{K_1}{T}$, $\frac{\widetilde F}{T}$ satisfy the conditions in Theorem \ref{CFSZ3}, so we have
\[
 t(C_5,K_1/T)\le t(C_5,\widetilde F/T)+11\cdot98^2(3\epsilon/T)^{1/4}
             =11\cdot98^2(3\epsilon/T)^{1/4}.
\]
Here $ t(C_5,\widetilde F/T)=0$ since the support of $F$ has no closed $5$-walk. When $\epsilon$ is sufficiently small, this is below the $\delta$ in
Theorem~\ref{CFSZ2} with $\zeta=\eta/(2T)$. The we can  find a symmetric $\mathcal D\subseteq V^2$ such that
\[
 t(C_5,(K_1/T)\mathbf1_{\mathcal D^c})=0,\qquad
 \frac1{n^2}\sum_{(x,y)\in\mathcal D}\frac{K_1(x,y)}{T}
 \le\frac{\eta}{2T}.
\]
Let $k=k_1\mathbf1_{\mathcal D^c}$. Since $\mathcal D$ is symmetric,
\[
\begin{aligned}
 K(x,y)&=k(x,y)+k(y,x)=K_1(x,y)\mathbf1_{\mathcal D^c}(x,y),\\
 t(C_5,K)&=T^5t(C_5,(K_1/T)\mathbf1_{\mathcal D^c})=0.
\end{aligned}
\]
This means the support of $k$ has no closed $5$-walk. The weight lost is
\[
 \frac1{n^2}\sum_{x,y}(k_1-k)(x,y)
 =\frac1{n^2}\sum_{(x,y)\in\mathcal D}k_1(x,y)
 =\frac1{2n^2}\sum_{(x,y)\in\mathcal D}K_1(x,y)\le\frac{\eta}{4}.
\]
We are done.
\end{proof}

Choose $T\ge\max\{1,16/\eta\}$. From above construction,  we have
$\mathcal{C}_k\le\mathcal{C}_{k_1}\le1$. By Lemma~\ref{lem:bookreduction},
\eqref{eq:bounded-weights} and Claim~\ref{cl:removing-walks},
\[
 \frac1{n^2}\sum_{x,y}k(x,y)
 \ge\frac{2}{3\sqrt{3}}+\eta-96\epsilon-\frac4T-\frac{\eta}{4}
 \ge\frac{2}{3\sqrt{3}}+\frac{\eta}{4}.
\]
The proof of this proposition is done.
\end{proof}

\section{The bipartite estimate}\label{sec:bipartite}
\begin{proposition}\label{lem:bipartite}
Let $V=X\sqcup Y$ be a set of size $n$. If $k\ge0$, $\mathcal{C}_k(x,y)\le 1$ and $k(x,y)=0$ when $(x,y)\in X^2\cup Y^2$, then
$n^{-2}\sum_{x,y}k(x,y)\le\frac{2}{3\sqrt{3}}$.
\end{proposition}
\begin{proof}
Suppose $p=|X|/n\ge q=|Y|/n$. If $q=0$, then $k(x,y)=0$ and there is nothing to prove. 
Let $t=n^{-2}\sum_{x,y}k(x,y)$. 
Expending the inequality $\sum_{y\in Y}\mathcal{C}_k(y,y)\le qn$, we have
\[
\begin{aligned}
 q&\ge\frac1{n^2}\sum_{x\in X,y\in Y}
\bigl(k(y,x)^2+\frac12k(x,y)^2+2k(x,y)k(y,x)\bigr)\\
 &\ge\frac1{2n^2}\sum_{x\in X,y\in Y}K(x,y)^2\ge\frac1{2n^2|X||Y|}
       \left(\sum_{x\in X,y\in Y}K(x,y)\right)^2
   =\frac{t^2}{2pq}.
\end{aligned}
\]
If $p\ge2/3$, then
$t^2\le2p(1-p)^2\le4/27$ and we are done. 

Therefore we may assume $q>0$ and  $1/2\le p<2/3$. Let $d:=2-3p$. So $0<d\le1/2$ and define the average outgoing and incoming weights
$u^+(x)=n^{-1}\sum_yk(x,y)$ and $u^-(x)=n^{-1}\sum_yk(y,x)$.
Define
\begin{equation*}
\begin{aligned}
 I:=&\frac4n\sum_{x\in X}\bigl(u^+(x)+(1-\frac{d}{8})u^-(x)\bigr)^2+\frac1n\left(1+\frac{7d}{4}+\frac{d^2}{64}\right)
                    \sum_{y\in Y}(u^+(y)+u^-(y))^2\\
 &+\frac{d(16-d)}{64n}\sum_{y\in Y}(2u^+(y)+u^-(y))^2
  +\frac{d(16-5d)}{64n}\sum_{y\in Y}u^+(y)^2.
\end{aligned}
\end{equation*}

\begin{Claim}
    $I\le (2+2d)p^2+4q^2.$
\end{Claim}
\begin{proof}
 We classify all the $2$-walks by their orientation. Call an arc type $F$ if it is from $X$ to $Y$ and type  $R$ otherwise.
Thus, there are four types of 2-walk: $FF,RR,FR,RF$. The four $2$-walks connecting two vertices in $X$ have coefficients list $r_X:=(1,1/2,1,1)$ as defined in Section 2. Analogously, the four $2$-walks connecting two vertices in $Y$ have coefficients list $r_Y:=(1/2,1,1,1)$.

For these four types of $2$-walks, define the coefficients
\[
\begin{aligned}
 s_X&=\left(3+2d,1-d+\frac{d^2}{8},3+\frac{7d}{4}+\frac{d^2}{64},
       3+\frac{7d}{4}+\frac{d^2}{64}\right),\\
 s_Y&=\left(0,4+d-\frac{d^2}{16},4+\frac{d}{2},4+\frac{d}{2}\right).
\end{aligned}
\]
Define the new $2$-walk capacity by
\[
\begin{aligned}
\text{for $a,b\in X$,} \quad S_k(a,b):=\frac1n\sum_{z\in Y}\Bigl[
&s_X(FF)k(a,z)k(b,z)
+s_X(RR)k(z,a)k(z,b)\\
&+s_X(FR)k(a,z)k(z,b)
+s_X(RF)k(z,a)k(b,z)
\Bigr],
\end{aligned}\]
\[
\begin{aligned}
\text{for $a,b\in Y$,} \quad S_k(a,b):=\frac1n\sum_{z\in X}\Bigl[
&s_Y(FF)k(z,a)k(z,b)
+s_Y(RR)k(a,z)k(b,z)\\
&+s_Y(FR)k(z,a)k(b,z)
+s_Y(RF)k(a,z)k(z,b)
\Bigr],
\end{aligned}
\]
and $S_k(a,b)=0$ when $a,b$ in different parts.  Here, for example $s_X(FF)$ is the coefficient in $s_X$ corresponding the 2-walk of type $FF$ .  Let the total capacity be $S=\frac1{n^2}\sum_{a,b\in V}S_k(a,b)$.

When consider a fixed vertex $y\in Y$, $\tfrac{1}{n^2}\sum_{x,x'\in X}k(x,y)k(x',y)=u^-(y)^2$,  similarly for other terms and then we get,
\[
\begin{aligned}
 \frac1{n^2}\sum_{x,x'\in X}\mathcal{C}_k(x,x')
 &=\frac1n\sum_{y\in Y}
       \bigl(u^-(y)^2+\tfrac12u^+(y)^2+2u^+(y)u^-(y)\bigr),\\
 \frac1{n^2}\sum_{y,y'\in Y}\mathcal{C}_k(y,y')
 &=\frac1n\sum_{x\in X}
       \bigl(\tfrac12u^+(x)^2+u^-(x)^2+2u^+(x)u^-(x)\bigr).
\end{aligned}
\]
The same discussion gives
\[
\begin{aligned}
 S=&\frac1n\sum_{x\in X}
       \bigl(s_Y(RR)u^-(x)^2+2s_Y(RF)u^+(x)u^-(x)\bigr)\\
 +&\frac1n\sum_{y\in Y}
       \bigl(s_X(FF)u^-(y)^2+s_X(RR)u^+(y)^2+2s_X(FR)u^+(y)u^-(y)\bigr).
\end{aligned}
\]
Expanding the squares defining $I$ and comparing these coefficients yields
\begin{equation}\label{eq:singlecount}
 S=\frac{4+4d}{n^2}\sum_{x,x'\in X}\mathcal{C}_k(x,x')
   +\frac8{n^2}\sum_{y,y'\in Y}\mathcal{C}_k(y,y')-I.
\end{equation}

Let $J:= \frac{1}{n^2}\sum_{a,b\in V}S_k(a,b)\mathcal{C}_k(a,b)$. We know $J\le S$ since each $\mathcal{C}_k(a,b)\le 1$ and $S_k(a,b)\ge 0$.
We  will compare $J$ with the  following $B_0$ and show $B_0\le J$:
\[
 B_0=\frac{2+2d}{n^2}\sum_{x,x'\in X}\mathcal{C}_k(x,x')^2
        +\frac4{n^2}\sum_{y,y'\in Y}\mathcal{C}_k(y,y')^2.
\]
Fix an ordered quadruple $(x,y,x',y')\in X\times Y\times X\times Y.$
Label its four edge positions by
\[
1=(x,y),\qquad 2=(x',y),\qquad
3=(x',y'),\qquad 4=(x,y').
\]
Each position has two possible directions, giving 16 direction choices. Let $E\subseteq\{1,2,3,4\}$ be the set of positions that the arc is type $F$ and $P_E(x,y,x',y')$ for the corresponding product of four arc weights. 
For a direction choice \((e_1,e_2,e_3,e_4)\), the term \(n^{-4}P_E(x,y,x',y')\) has coefficients
\[
\begin{aligned}
j_E
&=s_X(e_1e_2)r_X(e_4e_3)
  +s_Y(e_1e_4)r_Y(e_2e_3),\\
b_E
&=(2+2d)r_X(e_1e_2)r_X(e_4e_3)
  +4r_Y(e_1e_4)r_Y(e_2e_3)
\end{aligned}
\]
in $J$ and $B_0$, respectively. Thus
\[
\begin{aligned}
J&=\frac1{n^4}
\sum_{\substack{x,x'\in X\\y,y'\in Y}}
\sum_{E\subseteq\{1,2,3,4\}}
j_E P_E(x,y,x',y'),\\
B_0&=\frac1{n^4}
\sum_{\substack{x,x'\in X\\y,y'\in Y}}
\sum_{E\subseteq\{1,2,3,4\}}
b_E P_E(x,y,x',y').
\end{aligned}
\]
We now combine corresponding terms from different ordered quadruples. Define
\[
W_E:=\frac1{n^4}
\sum_{\substack{x,x'\in X\\y,y'\in Y}}
P_E(x,y,x',y')\ge0.
\]
Then
$J-B_0=\sum_{E\subseteq\{1,2,3,4\}}(j_E-b_E)W_E.$

For example,
\[
\begin{aligned}
P_{12}(x,y,x',y')
=k(x,y)k(x',y)k(y',x')k(y',x)=P_{34}(x,y',x',y).
\end{aligned}
\]
Since $(x,y,x',y')\longmapsto(x,y',x',y)$
is a bijection of the set of ordered quadruples being summed, we have \(W_{12}=W_{34}\). Hence the combined contribution of these two patterns is
\[
\begin{aligned}
(j_{12}-b_{12})W_{12}+(j_{34}-b_{34})W_{34}=\bigl(j_{12}-b_{12}+j_{34}-b_{34}\bigr)W_{12}.
\end{aligned}
\]
Substituting the values of \(r_X,r_Y,s_X,s_Y\), we obtain
\[
\begin{aligned}
j_{12}-b_{12}
&=\frac{3+2d}{2}+\left(4+\frac d2\right)-(5+d)
=\frac12+\frac d2,\\
j_{34}-b_{34}
&=\left(1-d+\frac{d^2}{8}\right)
  +\left(4+\frac d2\right)-(5+d)
=-\frac{3d}{2}+\frac{d^2}{8}.
\end{aligned}
\]
Their combined contribution is therefore
 $\left(\frac12-d+\frac{d^2}{8}\right)W_{12}\ge0.$
Similarly, exchanging $x,x'$ or exchanging $y,y'$ in the ordered quadruple $(x,y,x',y')$ gives
\[
\begin{gathered}
W_{123}=W_{124}=W_{134}=W_{234},\\
W_1=W_2=W_3=W_4,\\
W_{14}=W_{23},\qquad
W_{12}=W_{34},\qquad
W_{13}=W_{24}.
\end{gathered}
\]
Adding the coefficients of the patterns with equal weight sums gives the following table.

\begin{table}[htbp]
\centering
\setlength{\belowcaptionskip}{6pt}
\caption{Sums of the coefficient differences on the seven direction classes.}
\label{tab:pairs}
\small
\begin{tabular}{@{}ll@{}}
\toprule
F-position patterns & Sum of the coefficients $j_E-b_E$\\
\midrule
$1234$ & $0$\\
$\varnothing$ & $0$\\
$123,124,134,234$ & $d^2/32$\\
$1,2,3,4$ & $(1-d/8)(1-9d/8)$\\
$14,23$ & $0$\\
$12,34$ & $1/2-d+d^2/8$\\
$13,24$ & $2(1+d/8)^2$\\
\bottomrule
\end{tabular}
\end{table}

Consequently,
\[
\begin{aligned}
J-B_0
={}&\frac{d^2}{32}W_{123}
+\left(1-\frac d8\right)\left(1-\frac{9d}{8}\right)W_1\\
&+\left(\frac12-d+\frac{d^2}{8}\right)W_{12}
+2\left(1+\frac d8\right)^2W_{13}
\ge0.
\end{aligned}
\]

Using~\eqref{eq:singlecount} and $B_0\le J\le S$, we have 
\[\begin{aligned}
 I&\le\frac{2+2d}{n^2}\sum_{x,x'\in X}
       \bigl(2\mathcal{C}_k(x,x')-\mathcal{C}_k(x,x')^2\bigr)+\frac{4}{n^2}\sum_{y,y'\in Y}
       \bigl(2\mathcal{C}_k(y,y')-\mathcal{C}_k(y,y')^2\bigr)\\
 &\le \frac{2+2d}{n^2}|X|^2+\frac{4}{n^2}|Y|^2=(2+2d)p^2+4q^2.
\end{aligned}\]
This proves the claim.
\end{proof}

Next we  bound $I$ from below.
\begin{Claim}
  $I\ge\frac{27}{4}[(2+2d)p^2+4q^2]t^2$.
\end{Claim}
\begin{proof}

Let $\beta=n^{-2}\sum_{x\in X,y\in Y}k(y,x)$. So 
\[
 \frac1n\sum_{x\in X}u^+(x)=\frac1n\sum_{y\in Y}u^-(y)=t-\beta,
 \qquad
 \frac1n\sum_{x\in X}u^-(x)=\frac1n\sum_{y\in Y}u^+(y)=\beta.
\]
By above notation, we know the four terms in $I$ satisfy
\[
\begin{aligned}
&\frac{1}{n}\sum_{x\in X}
\left(u^+(x)+\left(1-\frac d8\right)u^-(x)\right)=t-\frac{d\beta}{8},\quad \frac{1}{n}\sum_{y\in Y}(u^+(y)+u^-(y))
=t, \\
&\frac{1}{n}\sum_{y\in Y}(2u^+(y)+u^-(y))
=t+\beta,\quad \frac{1}{n}\sum_{y\in Y}u^+(y)=\beta.
\end{aligned}
\]
Use Cauchy--Schwarz to each term, we get
\[
\begin{aligned}
 I
 &\ge\frac4p(t-d\beta/8)^2
  +\frac1q\left[
  \left(1+\frac{7d}{4}+\frac{d^2}{64}\right)t^2
  +\frac{d(16-d)}{64}(t+\beta)^2
  +\frac{d(16-5d)}{64}\beta^2\right]\\
 &=\frac{4}{p}(t-d\beta/8)^2
  +\frac{(1+2d)t^2+\frac{d(16-d)}{32}t\beta+\frac{d(16-3d)}{32}\beta^2}{q}.
\end{aligned}
\]
Let  $\Delta=\frac{4}{p}(t-d\beta/8)^2
  +\frac{(1+2d)t^2+\frac{d(16-d)}{32}t\beta+\frac{d(16-3d)}{32}\beta^2}{q}-\frac{27}{4}[(2+2d)p^2+4q^2]t^2$ and consider $12pq\Delta$.
Substitution of $p=(2-d)/3$ and $q=(1+d)/3$, the coefficients of $t^2,t\beta,\beta^2$ are, respectively,
$2d^3(d^2-2d+3)$, $-d^2(50-d)/8$ and
$d(5d^2-20d+32)/8$.
Hence 
\[
 12pq\Delta=d\left[
 2d^2(d^2-2d+3)t^2
 -\frac{d(50-d)}8t\beta
 +\frac{5d^2-20d+32}{8}\beta^2
 \right].
\]
For $0<d\le1/2$,
$2(d^2-2d+3)=2((1-d)^2+2)\ge9/2$, $(50-d)/8\le27/4$ and $\frac{5d^2-20d+32}{8}\ge \frac{81}{32}.$
Also since $t,\beta\ge0$, the expression in brackets is at least
\[
 \frac{9}{2}d^2t^2-\frac{27}{4}dt\beta+\frac{81}{32}\beta^2
 =\frac{9}{32}(4dt-3\beta)^2
 \ge0.
\]
This proves
$I\ge\frac{27}{4}[(2+2d)p^2+4q^2]t^2$.
\end{proof}

Now, combining above two claims,  we have  $\frac{27}{4}[(2+2d)p^2+4q^2]t^2\le I\le (2+2d)p^2+4q^2$, which implies $t\le\frac{2}{3\sqrt{3}}$. The proof is complete.
\end{proof}

\section{Completion of the proof}\label{sec:completion}
\begin{proof}[Proof of Theorem~\ref{thm:main}]
The lower bound is the construction of Bollob\'as and Gy\H ori~\cite{BG}.
Suppose that the upper bound fails.
By Proposition~\ref{prop:reduction}, there is a function $k$ on a finite set $V$ of size $n$ with
$t=n^{-2}\sum_{x,y}k(x,y)>\frac{2}{3\sqrt{3}}$, $\mathcal{C}_k\le1$, and no closed $5$-walk in its support. The support has no loops, since traversing a loop five times gives a closed $5$-walk.

Let $K(x,y)=k(x,y)+k(y,x)$ and $w=n^{-2}|\{(x,y)\in V^2:K(x,y)>0\}|$.
Applying Cauchy--Schwarz to the positive entries of $K$ and~\eqref{eq:capacitybounds}, we obtain
\[
 (2t)^2
 =\left(\frac1{n^2}\sum_{x,y}K(x,y)\right)^2
 \le w\|K\|_2^2\le\frac43w,
 \qquad w\ge3t^2>\frac49.
\]
Let $N(x)$ be neighborhood of $x$ in the support and $d(x)=|N(x)|$.
Then $n^{-2}\sum_xd(x)=w$.
Summing over ordered adjacent pairs yields
\[
 \frac1{n^3}\sum_{\substack{x,y\\K(x,y)>0}}(d(x)+d(y))
 =\frac2{n^3}\sum_xd(x)^2
 \ge2\left(\frac1{n^2}\sum_xd(x)\right)^2
 =2w^2.
\]
There are $wn^2$ ordered adjacent pairs, so some edge $uv$ satisfies
$d(u)+d(v)\ge2wn$.
If $z\in N(u)\cap N(v)$, then $u,v,z,u,v,u$ is a closed $5$-walk.
Thus $N(u)$ and $N(v)$ are disjoint. Consequently, for $B:=N(u)\cup N(v)$,
$|B|=d(u)+d(v)\ge2wn>8n/9$.

Define
$X=\bigcup_{a\in N(u)}N(a)$ and
$Y=\bigcup_{b\in N(v)}N(b)$.
If $x\in X\cap Y$, choose $a\in N(u)$ and $b\in N(v)$ adjacent to $x$, then
the sequence $u,a,x,b,v,u$ is a closed $5$-walk, a contradiction.
If two vertices $x,x'\in X$ are adjacent, choose $a,a'\in N(u)$ with $ax,a'x'$ being edges.
Then $u,a,x,x',a',u$ is also a closed $5$-walk.
The same argument holds for $Y$. Thus $X\cap Y=\varnothing$ and $X,Y$ are independent.

Let $U=X\cup Y=\bigcup_{b\in B}N(b)$, $\overline{U}=V\setminus U$ and $\rho=|\overline{U}|/n$. Since $u,v\in B$, we have $B=N(u)\cup N(v)\subseteq U$ and hence $0\le\rho\le1/9$.
By the definition of $U$, every $x\in\overline{U}$ has
$N(x)\subseteq V\setminus B$ and thus $|N(x)|\le n/9$.
Using Cauchy--Schwarz on $N(x)$ and then~\eqref{eq:capacitybounds} gives
\[
 \left(\frac1n\sum_yK(x,y)\right)^2
 \le\frac{|N(x)|}{n}\,\frac1n\sum_yK(x,y)^2
 \le\frac29
 \qquad(x\in \overline{U}).
\]
Thus the weight lost by deleting $\overline{U}$ satisfies
\begin{equation}\label{loss}
  t-\frac1{n^2}\sum_{x,y\in U}k(x,y)
 =\frac1{n^2}\sum_{x\in\overline{U},\,y\in V}K(x,y)
      -\frac1{n^2}\sum_{x,y\in\overline{U}}k(x,y)\le\frac{\sqrt2}{3}\rho.   
\end{equation}

On $U^2$, consider the new weight function $k'(x,y)=\sqrt{1-\rho}k(x,y)$.
Since $|U|=n(1-\rho)$, 
\[\mathcal{C}_{k'}(x,y)=\frac{1}{n}\sum_{z\in U}\left(k(x,z)k(y,z)+\tfrac12k(z,x)k(z,y)+k(x,z)k(z,y)+k(z,x)k(y,z)\right).\]
This implies $\mathcal{C}_{k'}(x,y)\le \mathcal{C}_{k}(x,y)\le 1$.
Also since $X,Y$ are independent, using Proposition~\ref{lem:bipartite} we get
\[
 \frac{\sqrt{1-\rho}}{|U|^2}\sum_{x,y\in U}k(x,y)\le\frac{2}{3\sqrt{3}},
 \qquad
 \frac1{n^2}\sum_{x,y\in U}k(x,y)\le\frac{2}{3\sqrt{3}}(1-\rho)^{3/2}.
\]
Combining this with the inequality (\ref{loss}), we have 
\[
 t\le\frac{2}{3\sqrt{3}}(1-\rho)^{3/2}+\frac{\sqrt2}{3}\rho
 \le\frac{2}{3\sqrt{3}}.
\]
The last inequality follows since the function  is devreasing when $0\le\rho\le1/9$.
This contradicts $t>\frac{2}{3\sqrt{3}}$ and the proof is completed.
\end{proof}

\section{Proofs of the applications}\label{sec:applications}
 In this section, we give short proof for the two applications.
\begin{proof}[Proof of Theorem~\ref{thm:induced}]
One may note, the same construction from Bollob\'as--Gy\H{o}ri construction is also induced-$C_4$-free, hence it gives the lower bound. For the proof of upper bound, only  Lemma \ref{deletion} need to be modified: during the process, we  delete the vertex of  degree less than $\frac{2}{3\sqrt{3}}(i^{3/2}-(i-1)^{3/2})$ and  retain all edges that lie in no triangle and orient them arbitrarily. The
deletion of $K_4$ blocks and book spines costs $o(n^{3/2})$ edges, so the
resulting oriented graph has $e(G)-o(n^{3/2})$ arcs and the similar discussion would show the Lemma \ref{lem:bookreduction} still hold for the directed graph.
Sections~\ref{sec:weighted}--\ref{sec:completion} then give the upper bound.\end{proof}

\begin{proof}[Proof of Theorem~\ref{thm:berge}]
Let $G$ be the shadow of a $\mathcal B$-free $3$-uniform hypergraph
$\mathcal H$, so a pair is an edge of $G$ if it is contained in a
hyperedge of $\mathcal H$. It is easy to show that $G$ must be $C_5$-free. Each hyperedge gives a distinct triangle in $G$, so
Theorem~\ref{thm:main} gives the upper bound. For the lower bound,
view each triangle  in the
Bollob\'as--Gy\H{o}ri construction   as  a hyperedges. This hypergraph is $\mathcal{B}$-free and hence give the lower bound.
\end{proof}

\section*{Declaration}
During the preparation of this work, the authors used AI systems for various tasks, including generation of candidate
proofs. At first the AI produced a very complicated proof, which the authors eventually streamlined into the current version. All AI
suggestions that were incorporated into the final version were verified and streamlined by the authors. The authors take full responsibility for the correctness of the paper.

\end{document}